\documentclass[a4paper,12pt]{amsart}

\usepackage[all]{xy}

\usepackage{t1enc}
\def\cQ{\mathcal{Q}}

\def\cF{\mathcal{F}}
\def\cH{\mathcal{H}}

\def\cN{\mathcal{N}}

\def\ce{\mathcal{E}}

\newcommand{\om}{\omega}

\newtheorem{theorem}{Theorem}[section]
\newtheorem{lemma}[theorem]{Lemma}
\newtheorem{proposition}[theorem]{Proposition}

\theoremstyle{remark}

\newtheorem*{definition*}{\rm\bf Definition}

\newcommand{\nn}[1]{(\ref{#1})}

\def\endrk{\hbox{$|\!\!|\!\!|\!\!|\!\!|\!\!|\!\!|$}}

\def\sideremark#1{\ifvmode\leavevmode\fi\vadjust{\vbox to0pt{\vss
 \hbox to 0pt{\hskip\hsize\hskip1em
 \vbox{\hsize3cm\tiny\raggedright\pretolerance10000
  \noindent #1\hfill}\hss}\vbox to8pt{\vfil}\vss}}}%
                        
\author{A. Rod Gover } \email{gover@math.auckland.ac.nz} \title{Q curvature
  prescription; forbidden functions and the GJMS null space}
\begin{document}

\address{ARG: Department of Mathematics\\
  The University of Auckland\\
  Private Bag 92019\\
  Auckland 1\\
  New Zealand} \email{gover@math.auckland.ac.nz}



\begin{abstract}
On an even conformal manifold $(M,c)$, such that the critical GJMS
operator has non-trivial kernel, we identify and discuss the role of a
finite dimensional vector space $\cN(\cQ)$ of functions determined by
the conformal structure. Using these we describe an infinite
dimensional class of functions that cannot be the Q-curvature $Q^g$ for
any $g\in c$.  If certain functions arise in $\cN(\cQ)$ then $Q^g$
cannot be constant for any $g\in c$.
\end{abstract}

\maketitle
\pagestyle{myheadings}
\markboth{Gover}{Q curvature and forbidden functions from $\cN(P)$}

\section{introduction}

Branson's Q-curvature $Q^g$ is a natural scalar Riemannian invariant
defined on even dimensional manifolds; it was introduced in the far
sighted works \cite{BrO,BrSeoul,Tomsharp} of Branson-\O rsted and
Branson.  Algorithms for obtaining $Q^g$, and explicit formulae in
low dimensions, may be found in \cite{GoPetCMP,GrH}. The problem of
conformally prescribing Q-curvature is that of determining, on a fixed
even dimensional conformal structure $(M^n,c)$, which functions may be
the Q-curvature $Q^g$ for some $g\in c$; in other words it is the
problem of characterising the range of the map
\begin{equation}\label{map}
Q:c\to C^\infty (M) \quad \mbox{ given by } \quad g\mapsto Q^g~.
\end{equation}
This has drawn considerable attention recently
(e.g.\ \cite{Baird,BaFaRe,Brendle,DR,MalStr}), and especially the uniformisation type
problem of whether one may find a metric in $c$ with $Q^g$ constant,
see for example \cite{CGY,CY,DM} and references therein.

Throughout we shall work on an even closed (that is, compact without
boundary) conformal manifold $(M,c)$. For simplicity of exposition we
shall assume that this is connected.  All structures will be taken
smooth and here $c$ is a Riemannian conformal class. That is it is an
equivalence class of Riemannian metrics where $g\sim \widehat{g}$
means that $\widehat{g}=e^{2\om}g$ for some $\om\in C^\infty (M)$. For
metrics related in this way, their Q-curvatures are related by
\begin{equation}\label{pQ}
Q^{\widehat{g}} = e^{-n \om} \big( Q^g+P^g \om \big),
\end{equation}
where the linear differential operator $P^g: C^\infty(M)\to C^\infty
(M)$ is the {\em critical} (meaning dimension order) {\em GJMS
  operator} of \cite{GJMS}. In dimension 4 the operator is due to
Paneitz, while in dimension 2 $P^g$ it is simply the Laplacian. In all
(even) dimensions this operator takes the form
$P^g=\Delta^{n/2}+{lower~order ~terms}$ where $\Delta$ indicates the
Laplacian. More important for our current purposes is that it is
formally self-adjoint \cite{GrZ} and may be expressed as a composition
\begin{equation}\label{Pform}
P^g=\delta M d,
\end{equation}
where $\delta$ the formal adjoint of the exterior derivative $d$, and
$M$ is some other linear differential operator (from 1-form fields to
1-form fields). Note also that $P^g$ is conformally invariant in the
sense that if $g$ and $\widehat{g}$ are conformally related, as above,
then $P^{\widehat{g}}=e^{-n\om}P^g$; in fact this is a necessary
consequence of the Q-curvature transformation law \nn{pQ}. Given a
linear operator $Op$, we shall generally write $\cN (Op)$ for its null
space.  From the conformal transformation law for $P^g$, it follows
that its null space $\cN(P^g)$  is conformally stable. That is
for any pair of metrics $g$, $\widehat{g}$ in the conformal class $c$
we have $\cN(P^g)=\cN(P^{\widehat{g}} )$. For this reason we shall
usually simply write $\cN(P)$ for this space.

 For the exterior derivative on functions, $\cN(d)$ is the space of
 constant functions.  Clearly \nn{Pform} implies $\cN(d)\subseteq
 \cN(P)$.  Much of the previous work on the prescription equation
 \nn{pQ} has assumed, or restricted to settings where, the operator
 $P$ has {\em trivial kernel}, meaning $\cN(d)=\cN(P^g)$.  In
 dimension $n=2$, $Q^g$ is (a constant multiple of ) the Gauss
 curvature and $P^d$ is the Laplacian. Thus in this case $P^g$ has
 trivial kernel on closed manifolds. However in higher dimensions the
 operator $P^g$ can have non-trivial kernel. For the Paneitz operator,
 in dimension 4, explicit examples are due to Eastwood and Singer
 \cite{ES-Frol}; see also \cite{ESlov-P}. (Each of these is a
 conformally flat product of a sphere with a Riemann surface and
 admits a metric $g$ with $Q^g$ constant but not zero.)  Thus here we
 restrict to $n\geq 4$ and make some first steps to study the
 consequences of the GJMS operator $P^g$ having non-trivial kernel.
 The first such is the existence of a class of forbidden functions,
 that is functions that cannot be in the range of $Q$ (in \nn{map}),
 as follows.
\begin{theorem}\label{forbidden}
  On a closed connected conformal manifold $(M,c)$ with $\cN(P)\neq
  \cN(d)$, there are non-zero functions $u\in \cN(P)$ such that for no
  pair $(g\in c,\alpha\in \mathbb{R})$, $\alpha \neq 0$, we have $Q^g
  =\alpha u$. If $\int_M Q$ is not zero or if $\dim \cN(P)\geq 3$ then
  there are such functions $u$ which are non-constant.
 In any case there is
  a space of functions $\cF$, with infinite dimensional linear span,
  such that if $ f\in \cF$ then for no pair $(g\in c,\alpha\in
  \mathbb{R})$, $\alpha \neq 0$, we have $Q^g =\alpha f$.
\end{theorem}

Most of the Theorem is established in Section \ref{calcs}. Lemma
\ref{noQ} identifies a conformally determined linear space of
functions (denoted $\cN(\cQ)$) as forbidden functions in
$\cN(P)$. These functions yield more general constraints on the range
of $Q$, as described in Propositions \ref{noQ2} and \ref{strict}.
Section \ref{role} explores the space $\cN(\cQ)$ and its
relation to $\cN(P)$. Proposition \ref{cap} shows that the span of
this subspace of $\cN(P)$ is a certain ``orthogonal complement'' to
the range of $Q$ in \nn{map}.  
Proposition \ref{idforbid} gives a strengthening of Theorem
\ref{forbidden} for the cases where $\int Q$ is non-zero.  Theorem
\ref{decthm} shows that, when $\int Q$ is not zero, $\cN(P)$ admits a
conformally invariant direct sum decomposition into the space of
constant functions versus the functions suitably orthogonal to $Q$. As
explained there, one reason this is interesting is because, although
this decomposition is conformally invariant, in the case that there is
a metric $g$ with $Q^g$ constant it recovers the Hodge decomposition
of $\cN(P)$.  
The main remaining result in Section \ref{role} is Theorem
\ref{1notin}. In the cases $\int Q\neq 0$ this shows that certain
functions in $\cN(\cQ)$ (if they arise) obstruct the prescription of
constant Q-curvature. This and the (more obvious) analogous result for
$\int Q=0$ lead to two interesting open questions which are posed.

In Section \ref{harm} we explain the link between the results here and
the so-called conformal harmonics of \cite{BrGodeRham}.  Finally in
Section \ref{linQ} we point out that almost all results in the article
will extend to the prescription problem for natural invariants in the
so-called linear Q-space. The latter are invariants which 
transform conformally suitably like the Q-curvature.

On the sphere forbidden functions for Q-prescription have been
identified by Delano\"{e} and Robert \cite{DR} (and in fact those
authors consider a broader class of invariants than what we here call
Q-curvature).  These functions violate symmetry related 
Kazdan-Warner type identities necessarily satisfied by $Q^g$; these
identities generalise those known for the scalar curvature
\cite{KW,B}. The functions concerned are obviously not in $\cN(P)$
and this is essentially a different phenomenon to that discussed
here. 

Some issues related to studying Q-prescription when $P^g$ has
non-trivial kernel were touched on in \cite{MalchSIGMA}. The current
work was partly motivated by the questions suggested there and by
discussions with the author of that work, Andrea
Malchiodi. Conversations with Michael Eastwood are also much
appreciated.

\subsection*{Acknowledgements} The basic idea for this work arose
during the conference ``PDE In Conformal and K\"ahler Geometry'',
University of Science \& Technology of China, July 2008 and also the
New Zealand Institute of Mathematics and Its Applications (NZIMA)
thematic programme ``Conformal Geometry and its Applications'' hosted
by the Department of Mathematics of the University of Auckland in
2008. The article was written during the programme ``Geometric Partial
Differential Equations'' at the Institute for Advanced Study,
Princeton. The author is grateful for the stimulation and support of
these programmes and institutions.  ARG is also supported by Marsden
Grant no.\ 06-UOA-029

\section{The Proof of Theorem \ref{forbidden}} \label{calcs}

We recall from \cite{BGaim} (updated 
as \cite{BGact}) the following observation.
\begin{proposition}\label{Qu}
For $u\in \cN(P)$,
$$
\int_M u Q^g  ~\mu_g
$$
is conformally invariant.
\end{proposition}
\noindent Here we write $\mu_g$ for the canonical measure determined
by the metric. For our current purposes it is useful to review the
proof of the Proposition.\\
\noindent{\bf Proof:} Recall that if $\widehat{g}=e^{2\omega}g$,
$\om\in C^\infty (M)$, then we have the Q-curvature transformation
\nn{pQ}.
 Since the measure $\mu_g$
transforms conformally according to
$$
\mu_{\widehat{g}}= e^{n\om} \mu_g
$$
we have 
$$
\int_M uQ^{\widehat{g}}~\mu_{\widehat{g}} = \int_M u e^{-n \om} \big(
Q^g+P^g \om \big) e^{n\om} ~\mu_g = \int_M u 
Q^g ~\mu_g + \int_M u (P^g \om)  ~\mu_g ~.
$$
But $P^g$ is formally self-adjoint \cite{GrZ} and so the second term drops out.
\quad $\Box$

Let us write $\cQ$ for the linear operator $\cQ:\cN(P)\to \mathbb{R}$ given by 
$$
u\mapsto \int_M Q^g u ~\mu_g,
$$ 
and $\cN(\cQ)$ for its kernel.
Proposition \ref{Qu} may be viewed as the statement that $\cQ$ is
conformally invariant.

\begin{lemma}\label{noQ} If \hspace{1pt} $0\neq u\in \cN(\cQ)$ then there is no
  $g\in c$ with $Q^g=u$.
\end{lemma}
\noindent{\bf Proof:} Given $u\in \cN(\cQ)$, suppose that for some $g'
\in c$ we have $Q^{g'}=u $. Then $u \in \cN(\cQ) $
implies that
$$
0=\int_M u Q^g  \mu_g .
$$ Since $u\in \cN(P)$, the right hand side is conformally
invariant. So we may calculate in the metric $g'$, whence
$$
0= \int_M u Q^{g'}  \mu_{g'}= \int_M |u|^2 \mu_{g'}
$$
and so $u=0$.
\quad $\Box$

\noindent{\bf Proof of the first statement in Theorem
  \ref{forbidden}:} Since $\cN(d)\subseteq \cN(P)$, and $\cQ:\cN(P)\to
\mathbb{R}$ is a linear map it follows that $\cN(\cQ)$ has codimension
at most one. Thus if $\cN(d)\neq \cN(P)$ then $\cN(\cQ)\neq \{
0\}$. So the first statement follows from the previous Lemma.
\quad $\Box$

So the non-zero elements of $\cN(\cQ)$ are forbidden functions. In
fact we should really view $\cN(\cQ)\setminus \{0\}$ (and more generally 
$\cN(P)\setminus \{0\}$) as a space of
constraints on the range of $Q$ as follows. 
If $f\in C^\infty (M)$ is in the
range of $Q$ then it must be that there is $g\in c$ with $\int_M f u
\mu_g=0$ for all $u\in \cN(\cQ)\setminus \{0\}$.  For $g\in c$, let us
write $\cN(\cQ)^{\perp_g}$ for the orthogonal complement in
$C^\infty(M)$ of $\cN(\cQ)$ with respect to the usual $L^2$ inner
product on $(M,g)$.  Let us make the definitions
$$
\cF^g:= \{ f\in C^\infty (M)~:~ f \notin \cN(\cQ)^{\perp_g}\}
$$
and 
$$
\cF:=\cap_{g\in c} \cF^g.
$$ 
From the definition of $\cF$ and an obvious adaption of the proof
of Lemma \ref{noQ} we have the following.
\begin{proposition}\label{noQ2} If $ f\in \cF$ then 
for no
  pair $(g\in c,\alpha\in \mathbb{R})$, $\alpha \neq 0$, we have $Q^g
  =\alpha f$.
\end{proposition}
\noindent From the proof of Lemma \ref{noQ} we  also have that 
\begin{equation}\label{Fbig}
\cF\supseteq \cN(\cQ)\setminus \{0\} .
\end{equation}

\medskip

The last statement of Theorem \ref{forbidden} uses the following.
\begin{proposition}\label{strict} If $\cN(\cQ)$ contains non-constant 
functions  then 
the containment in \nn{Fbig} is proper and the span of $\cF$ is infinite 
dimensional.
\end{proposition}
\noindent{\bf Proof:} Observe that if $u\in \cN(\cQ)$ then, for
example, $u^p\in \cF$ for $p$ an odd positive integer.  
If $u$ is non-constant then there is no
  linear relation among the  $u^p$ as $p$ ranges over odd positive integers.
On the other hand $\cN(\cQ)\subseteq \cN(P)$ and $\cN(P)$ is finite dimensional 
since $P$ is elliptic and $M$ is compact.
\quad $\Box$\\
\noindent In the spirit of the proof here, note that if
$f:\mathbb{R}\to \mathbb{R}$ is any function with the same strict sign
as the identity then, for any $u\in \cN(\cQ)\setminus
\{0\}$, $f(u)$ is in $\cF$.

So to complete the proof of the Theorem, the main remaining task 
is to show when
$\cN(\cQ)$ necessarily contains non-constant functions.  Although
various cases are easily settled, it is worthwhile to first look at the
structure of $\cN(P)$.

\section{The structure and role  of $\cN(P)$ and $\cN(\cQ)$}\label{role}

Given a metric $g$ on $M$, let us write
$\overline{\cQ}^g:C^\infty(M)\to \mathbb{R}$ for the map 
$$
f\mapsto \int_M f Q^g~\mu_g.
$$ (So $\cQ$ agrees with the restriction of $\overline{\cQ}^g$ to
$\cN(P)$.)  
We have the
following interpretation of $\cN(\cQ)$.
\begin{proposition}\label{cap} Given a conformal class $c$ on a closed manifold $M$ we have
\begin{equation}\label{capQ}
\cN(\cQ)= \cap_{g\in c}\cN( \overline{\cQ}^g).
\end{equation}
\end{proposition}
\noindent{\bf Proof:} $\subseteq$: From the definitions of
$\overline{\cQ}^g$ and $\cQ$, and the conformal invariance of the
latter, it is immediate that for every $g\in c$ we have $\cN(\cQ)\subseteq \cN(
\overline{\cQ}^g)$.

\medskip

\noindent $\supseteq $:
Since $P^g$ is formally self-adjoint, given $u\in C^\infty (M)$, it
follows easily from \nn{pQ} that if
$$
\int_M u Q^g ~\mu_g=\int_M u Q^{\widehat{g}} ~\mu_{\widehat{g}} 
\quad \mbox{for all}\quad (g,\widehat{g})\in c\times c, 
$$ then $u\in \cN(P)$. If $u\in \cap_{g\in c}\cN( \overline{\cQ}^g)$
then we have the situation of the display, and also that $\int_M u Q^g 
\mu_g =0 $. So $u\in \cN(\cQ)$.  \quad $\Box$\\

\noindent{\bf Remark:} Note that, since the $L^2$ inner product is
definite, (excepting 0) the right hand side of \nn{capQ} consists of
functions which are manifestly not in the range of $Q$. Thus the
Proposition above gives an alternative proof of Lemma \ref{noQ}.

Here we are regarding $Q$ as function valued. However we might also
consider $Q$ as taking values in conformal densities of weight $-n$,
as in \cite{BrGodeRham}. As discussed there, there is a conformally
invariant pairing $\langle~,~\rangle$ between functions and such
densities. In terms of this 
$\langle \cN(\cQ),
~\rangle $ is the annihilator of the range of the map \nn{map} and
this characterises the space $\cN(\cQ)$.   \quad \endrk

For $q \in {\mathbb R}$ we use, as usual, the same notation for the constant
function on $M$ (which, recall, we assume connected) with value
$q$. Since, $P_k 1=0$, it follows from Proposition \ref{Qu} that, in
particular,
$$
k_Q:=\int_M Q^g~\mu_g
$$ is a global conformal invariant. Of course this was known well
before Proposition \ref{Qu} from \cite{BrO,BrSeoul}. This gives
immediate restrictions on the range of $Q^g$. Let us write
$\mathcal{E}_+$ (resp.\ $\mathcal{E}_-$) for the space of functions
$f\in C^\infty(M)$ such that $f$ is pointwise non-negative
(respectively non-positive) but not identically zero. We will write $\mathcal{E}_\pm$ for the union of these spaces.
 If $k_Q=0$ then we have the well known
result that any $f\in \mathcal{E_\pm}$
is not in the range of $Q$. In fact this is an example of $f\in
\cF\setminus \big( \cN(\cQ)\setminus \{0\} \big)$: when $k_Q=0$ we
have $1\in \cN(\cQ)$ and for no metric $g\in c$ is $f$ orthogonal to
1. Let us record then that in this case, without assuming the
containment $\cN(P)\supseteq\cN(d)$ is proper, we have the following.
\begin{proposition} \label{Qzero}
If $k_Q=0$ then $\cF\supseteq \mathcal{E_\pm}$. In particular $\cF$ 
spans an infinite dimensional vector space.
\end{proposition}
\noindent Thus the last statement of the Theorem \ref{forbidden} has no new
information in the case of $k_Q=0$.

When $k_Q\neq 0$ then we still obtain an obvious (and well known)
constraint from $1\in \cN(P)$, but this of a slightly different nature
since the non-zero constants are not in $\cN(\cQ)$.  Combining this
with the observations of Section \ref{calcs} we have the following.
\begin{proposition} \label{idforbid}
If $k_Q>0$ (alternatively $k_Q<0$) then the functions in
$\mathcal{E_-}\cup \cF\cup \{ 0\}$ (resp.\ $\mathcal{E_+}\cup \cF\cup
\{0 \}$) are not in the range of $Q$.
\end{proposition}

\medskip

Next we exhibit a decomposition of $\cN(P)$ which establishes the
second statement of Theorem \ref{forbidden} for the case that $k_Q\neq
0$; it shows that in this case the forbidden functions of $\cN(\cQ)$ are
necessarily non-constant. But it gives more than this and is of
independent interest.
\begin{theorem}\label{decthm}
  Suppose that $k_Q\neq 0$. Then we have a
  conformally invariant direct decomposition
\begin{equation}\label{dec}
\cN(P)= \cN(d)\oplus \cN\big( \cQ)~.
\end{equation}
\end{theorem}
\noindent{\bf Warning:} The decomposition of the Theorem is not necessarily orthogonal for any metric $g\in c$.\\
\noindent{\bf Proof of the Theorem:} $\cN(P)$ and $\cN(d) $ are conformally invariant
linear subspaces of $C^\infty (M)$. The vector space inclusion
$\cN(d)\hookrightarrow \cN(P) $ is split by the conformally invariant map
$$
\cN(P)\ni u\mapsto u_0:= \frac{1}{k_Q}\cQ(u).
$$
(So explicitly \nn{dec} is given by  $u= u_0+(u-u_0)$ for $u\in \cN(P)$.)
\quad $\Box$

\medskip

\noindent{\bf Remark:} Note that, for $\cN(P)$, \nn{dec} is a conformal
version of the Hodge decomposition. We mean this
as follows.

Suppose that there is $g\in c$ with $Q^g=q $ constant. Since we assume
$k_Q\neq 0$, it follows that $q\neq 0$.  By the Hodge decomposition on
$(M,g)$, for any $u \in C^\infty (M)$ we have
\begin{equation}\label{hodge}
u=\overline{u}+u_1
\end{equation}
where $\overline{u}$ is a constant function and $u_1$ is a divergence. In 
particular we may apply this decomposition to $u\in \cN(P)$. We have
$\overline{u}\in \cN(d)$, $u_1\in \cN(P)$ and 
$$
\int_M u_1 Q^g  \mu_g= q\int_M  u_1 \mu_g=0.
$$
So $u_1\in \cN(\cQ)$. On the other hand 
$$ u_0:=\frac{1}{k_Q}\int_M u Q^g ~\mu_g= \frac{q}{k_Q} \int_M
(\overline{u}+u_1 ) ~\mu_g =\frac{\overline{u}}{k_Q}\int_M q
~\mu_g =\overline{u}.
$$ So also we have $u_1=u-u_0$.  

This shows that on $(M,g)$ the Hodge decomposition \nn{hodge}, of
$\cN(P)$, agrees with \nn{dec}. But the latter is conformally
invariant and so gives a conformally invariant and canonical route to
the Hodge decomposition of $\cN(P)$ with respect to the metric $g$
that has $Q^g$ constant.

In general we do not know there is a metric $g$ that makes $Q^g$
constant; there is no preferred metric to exploit for a Hodge
decomposition. Nevertheless we always have the conformal decomposition
\nn{dec}. \quad \endrk

\medskip

These observations suggest the following problem. Recall $\ce_{\pm}$ is
the space of non-zero functions which are either non-negative or non-positive.

\smallskip

\noindent{\bf Question 1:} Suppose that $k_Q\neq 0$. Can the finite
dimensional vector space $\cN(\cQ)$ intersect non-trivially with
$\mathcal{E}_\pm$?

\smallskip

\noindent This question is interesting because if there are such
functions then they obstruct the prescription of constant Q-curvature.
To simplify the statement, note that $f\in \cN(\cQ)\cap \mathcal{E_-}$
is equivalent to $-f \in\cN(\cQ)\cap \mathcal{E_+} $.
\begin{theorem} \label{1notin}
Assume $(M,c)$ is a closed conformal manifold with $k_Q\neq 0$. 
Suppose  there is $u\in \cN(\cQ)\setminus \{ 0 \}$ such that $u$ is
pointwise non-negative. Then  $\not\exists g\in c$ with $Q^g$ either 
pointwise positive or pointwise negative. In particular there is no $g\in c$ with 
$Q^g$ constant.
\end{theorem}

\noindent{\bf Proof:} From Proposition \ref{Qu} $0=\int_M u Q^g \mu_g
$ is conformally invariant, but for no metric is $u$ is orthogonal to
to a strictly positive or strictly negative function.  \quad $\Box$

Note that if $k_Q=0$ then $ \cN(d)\subseteq \cN\big( \cQ).  $ So we
cannot hope to have the decomposition \nn{dec}. 
On the other hand in this case there is the
possibility that $\cN(\cQ)=\cN(P)$.  There is a characterisation of
this situation, as follows.
\begin{proposition}\label{Ddo} On a conformal manifold $(M,c)$ 
$$
\big(~ \cN(P)= \cN(\cQ)~ \big) \Leftrightarrow \big(~ \exists g\in c~ \mbox{s.t.}~ Q^g=0 ~\big)
$$
\end{proposition}

\noindent{\bf Proof:} $\Leftarrow :$ For all $u\in \cN(P)$, 
since we may calculate the conformal invariant $\cQ(u)=\int_M u Q^g  \mu_g$
using the metric $g\in c$  satisfying $Q^g=0$, it is clear that $\cQ(u)=0$.

\medskip

\noindent $\Rightarrow :$ Suppose we start from an
arbitrary metric $g\in c$ and want to find $\widehat{g}\in c$ with
$Q^{\widehat{g}}=0$. The prescription equation \nn{pQ} has
the simple form
$$
Q^g+P^g \om =0
$$ and so, by standard Fredholm theory, is solvable if and only if
$Q^g$ is orthogonal $\cN(P)$ i.e.\ if and only if $\cN(P)\subseteq \cN(\cQ) $.
\quad $\Box$

\noindent Note that Proposition \ref{Ddo} is simply a restatement of
an observation of Malchiodi in \cite[Section 3]{MalchSIGMA} and there the
following question is suggested.

\smallskip

\noindent{\bf Question 2:} Suppose that $k_Q=0$. Can there be
functions in $\cN(P)$ which are not orthogonal to $Q^g$?

\smallskip

\noindent As for question 1, such functions obstruct the prescription of constant $Q^g$. 
 
\noindent{\bf Remark and the proof of Theorem \ref{forbidden}:}
Suppose that $k_Q=0$ and $\dim \cN(P)=2$. Then either
$\cN(\cQ)=\cN(P)$, and $0$ is in the range of $Q$ but all other
elements of $\cN(P)$ are forbidden, or $\cN(\cQ)=\cN(d)$ in which case
we cannot solve for constant $Q^g$, but we cannot identify
non-constant forbidden functions in $\cN(P)$.  On the other hand if
$\dim \cN(P)> 2$ then it is clear that always we get non-constant (and
even mixed sign) functions in $\cN(P)$ that are not in the range of
$Q$ (as claimed in Theorem \ref{forbidden}). On the other hand when
$k_Q\neq 0$ it is immediate from, for example, Theorem \ref{decthm} that
$\cN(\cQ)\setminus \{0\}$ consists of non-contant functions.  \quad
\endrk

\section{Other links} \label{links}

We sketch here links with some related directions.

\subsection{Conformal Harmonics} \label{harm}

As above the setting is an even conformal manifold $(M^n,c)$.  In
\cite{BrGodeRham} a space of so-called {\em conformal harmonics}
$\cH^k$ is defined for each $k =1,\cdots n/2$ (see also \cite{E-Guill}). This
 is a conformally
stable subspace of $k$-forms that fits into an exact complex linking
it to the $k^{\rm th}$ de Rham cohomology space $H^k$
\cite[Proposition 2.5]{BrGodeRham}. For $k=1$, and in our current
notation, the complex is
\begin{equation}\label{cx}
0\to \cN(d)\to \cN(P) \stackrel{d}{\to} \cH^1 \to H^1
\end{equation}
where the map $\cN(P)\to \cH^1$ is just the restriction of $d$ and
$\cH^1\to H^1$ takes each 1-form in $\cH^1$ to its class in $H^1$.  It is
not known if the last map is necessarily surjective; by \cite[Theorem
  2.6]{BrGodeRham} it is if $\cN(d)= \cN(P)$ and then $\cH^1\cong
H^1$ (this is termed strong 0-regularity). 

Evidently the kernel of the map $\cH^1\to H^1$ measures the non-triviality of
the null space of the critical GJMS operator $P$. If $k_Q\neq 0$ then Theorem
\ref{decthm} shows that the complex in the display may be simplified
to
$$
0\to \cN(\cQ) \stackrel{d}{\to} \cH^1 \to H^1.
$$

As mentioned above, there is an analogue of the sequence \nn{cx} for
for each $k=0,1,\cdots ,n/2$. It turns out the ideas of section
\ref{role} can be adapted to give a generalisation of Theorem
\ref{decthm} which applies to all of these (and hence yields
additional structure to the theory in \cite{BrGodeRham}), although for
$k\geq 2$ the situation is rather more subtle than the case here.
This will be taken up elsewhere.

\subsection{Prescription problems for natural scalar invariants in the 
linear Q-space} \label{linQ}
Constructions of natural scalar Riemannian invariants
with conformal transformation properties similar to the Q-curvature
have been described in \cite{FeffH,GoPetCMP} and \cite{BrGoPont}. In
 Section 5 of the last reference two systematic
constructions are given for such invariants $K^g$ and in each case
this has the property that for $g,\widehat{g}\in c$ we have (cf. \nn{pQ})
$$ \widehat{g}=e^{2\om}g \quad \Rightarrow \quad
K^{\widehat{g}}=e^{-n\om}(K^g+L^g \om),
$$ $\om\in C^\infty (M)$, where $L^g$ is a formally self-adjoint
linear differential operator (necessarily conformally invariant)
of the form a composition
$$
L^g= \delta \tilde{M} d ,
$$ where, recall, $\delta$ the formal adjoint of $d$ and $\tilde{M}$
is some linear differential operator.  Each such quantity $K^g$ yields
a variational prescription problem (the case of prescribing constant
$K^g$ is discussed in Section 7 of \cite{BrGoPont}) along the same
lines as that for the $Q$-curvature.

It follows easily from these conformal transformation properties
displayed that Theorem \ref{forbidden} and, apart from Proposition
\ref{Ddo},  essentially all the results from Section \ref{calcs} above will
hold if $Q^g$ and $P^g$ are replaced by, respectively $K^g$ and
$L^g$. For Proposition \ref{Ddo} we still have the implication
$\Leftarrow$ if these replacements are made. The other adjustments are:
for the analogue of Proposition \ref{strict} we do not necessarily
 have proper containment and this would affect the analogue of Theorem
 \ref{forbidden} in the obvious way; since $L^g$ may not be elliptic, the space
$\cN(L^g)$ will not necessarily be finite dimensional.

Note that at one extreme we have the case that $K^g$ is simply a
natural conformal invariant (e.g.\ the Weyl curvature squared
$||W||^2$ in dimension 4). Then
$L^g$ is the zero
operator.  In this
case we trivially have $\cN(L)=C^\infty(M)$ and the analogue of
$\cN(\cQ)$ is strictly smaller than $\cN(L)$ and consists of all
functions which are orthogonal (with respect to the $L^2$ inner
product for $g$) to $K^g$ for all $g\in c$.


\begin{thebibliography}{XX}

\bibitem{E-Guill} E.\ Aubry, and C.\ Guillarmou, {\em Conformal
  harmonic forms, Branson-Gover operators and Dirichlet problem at
  infinity}. Preprint arXiv:0808.0552


\bibitem{Baird} P.\ Baird, A.\ Fardoun, and R.\  Regbaoui, {\em $Q$-curvature
  flow on 4-manifolds},  Calc.\ Var.\ Partial Differential Equations, {\bf 27}
  (2006),  75--104.

\bibitem{BaFaRe} P.\ Baird, A.\ Fardoun, and R.\ Regbaoui, {\em
  Prescribed Q-curvature on manifolds of even dimension}, preprint.

\bibitem{B} J.-P.\ Bourguignon, and J.-P. Ezin, {\em Scalar curvature
  functions in a conformal class of metrics and conformal
  transformations}, Trans.\ Amer.\ Math.\ Soc.\ {\bf 301} (1987), 
  723--736.

\bibitem{BrSeoul} T.P.\ Branson, The functional
  determinant, Lecture Notes Series, 4. Seoul National University,
  Research Institute of Mathematics, Global Analysis Research Center,
  Seoul, 1993, vi+103 pp.

\bibitem{Tomsharp} T.P.\ Branson, {\em  Sharp inequalities, the
  functional determinant, and the complementary series},
  Trans.\ Amer.\ Math.\ Soc.\  {\bf 347} (1995), 3671--3742.

\bibitem{BGaim} T.P.\ Branson, and A.R.\  Gover, {\em Origins,
  applications and generalisations of the $Q$-curvature}, American
  Institute of Mathematics, (2003), 
  http://www.aimath.org/pastworkshops/confstruct.html


\bibitem{BGact} T.P.\ Branson, and A.R. Gover, {\em Origins, applications
  and generalisations of the $Q$-curvature},  Acta Appl. Math. {\bf 102}
  (2008), 131--146.

\bibitem{BrO}
T.P.\ Branson, and B.\ \O rsted, 
{\em Explicit functional determinants in four dimensions},
Proc.\ Amer.\ Math.\ Soc., {\bf 113} (1991), 669--682. 


\bibitem{BrGodeRham} T.P.\ Branson, and A.R.\ Gover, {\em  Conformally
  invariant operators, differential forms, cohomology and a
  generalisation of $Q$-curvature},  Comm.\ Partial Differential
  Equations {\bf 30} (2005), 1611--1669.

\bibitem{BrGoPont} T.\ Branson, and A.R.\ Gover, {\em Pontrjagin forms and
  invariant objects related to the $Q$-curvature},
  Commun.\ Contemp.\ Math.\ {\bf 9} (2007),  335--358.

\bibitem{Brendle} S.\ Brendle, {\em Convergence of the $Q$-curvature flow
  on $S\sp 4$},  Adv.\ Math.\  {\bf 205} (2006), 1--32.


\bibitem{CGY} S.-Y.A.\ Chang,  M.J.\ Gursky, and P.C.\  Yang, {\em A 
conformally invariant sphere theorem in four dimensions},  Publ.\ Math.\ Inst.\
 Hautes \'Etudes Sci.\  No.\ 98  (2003), 105--143. 


\bibitem{CY} S.-Y.A.\ Chang, and P.C.\ Yang, {\em Extremal metrics of
  zeta function determinants on $4$-manifolds}, Ann.\ of Math.\ (2)
  {\bf 142} (1995), 171--212.


\bibitem{DR}
P.\ Delano\"e, and F.\ Robert,
{\em On the local Nirenberg problem for the $Q$-curvatures},
Pacific J.\ Math., {\bf 231} (2007), 293--304. 

\bibitem{DM}
Z.\ Djadli, and A.\ Malchiodi,
{\em Existence of conformal metrics with constant $Q$-curvature}, 
Ann.\ of Math. (2), to appear. \quad  
Preprint arXiv:math/0410141 


\bibitem{ES-Frol} M.G.\ Eastwood, and M.A.\ Singer, {\em The
  Fr\"ohlicher  spectral sequence on a twistor space},
  J.\ Differential Geom.\  {\bf 38} (1993), 653--669.

\bibitem{ESlov-P} M.G. Eastwood, and J.\ Slovak, {\em A primer on Q-curvature}, 
American
  Institute of Mathematics, (2003), \\
  http://www.aimath.org/pastworkshops/confstruct.html

\bibitem{FeffH} C.\ Fefferman, and K.\ Hirachi, {\em Ambient metric
  construction of $Q$-curvature in conformal and CR geometries}
  Math.\ Res.\ Lett.\ {\bf 10} (2003), 819--831.

\bibitem{GoPetCMP} A.R.\ Gover, and L.J.\ Peterson, {\em Conformally
  invariant powers of the Laplacian, $Q$-curvature, and tractor
  calculus},  Comm.\ Math.\ Phys.\  {\bf 235} (2003),  339--378.

\bibitem{GJMS} C.R.\ Graham, R.\ Jenne, L.J.\ Mason, and G.A.\ Sparling,
  {\em Conformally invariant powers of the Laplacian
  I. Existence}  J.\ London Math.\ Soc.\ (2) {\bf 46} (1992),
 557--565.

\bibitem{GrH} C.R.\ Graham, and K.\ Hirachi, {\em The ambient
  obstruction tensor and $Q$-curvature},  AdS/CFT correspondence:
  Einstein metrics and their conformal boundaries, 59--71, IRMA
  Lect. Math. Theor. Phys., 8, Eur. Math. Soc., Zürich, 2005.

\bibitem{GrZ} C.R.\ Graham, and M.\ Zworski, {\em Scattering matrix in
  conformal geometry},  Invent.\ Math.\  {\bf 152} (2003), 89--118. 

\bibitem{KW} J.L.\ Kazdan, and F.W.\  Warner, {\em  Scalar curvature and
  conformal deformation of Riemannian structure},  J. Differential
  Geometry {\bf 10} (1975), 113--134.

\bibitem{MalchSIGMA} A.\ Malchiodi, {\em Conformal metrics with constant
  $Q$-curvature},  SIGMA Symmetry Integrability Geom. Methods Appl.  {\bf 3}
  (2007), Paper 120, 11 pp.

\bibitem{MalStr} A.\ Malchiodi, and M.\ Struwe, {\em $Q$-curvature flow on $S\sp 4$},
J.\ Differential Geom.\ {\bf 73} (2006),  1--44. 

\end{thebibliography}
\end{document}